\documentclass[11pt]{article}
\usepackage[latin1]{inputenc}
\usepackage[T1]{fontenc}
\usepackage[english]{babel}
\usepackage{amsmath}
\usepackage{amsthm}
\usepackage{amsfonts}
\usepackage{amssymb}
\usepackage{graphicx}
\usepackage{pgf,tikz}
\usetikzlibrary{arrows}
\usepackage{stmaryrd}
\usepackage{tikz}
\usepackage{url}
\usepackage[np]{numprint}
\usepackage{vmargin}
\usepackage{mathrsfs}
\newtheorem{thm}{Theorem}
\newtheorem{lem}{Lemma}
\newtheorem{prop}{Proposition}

\theoremstyle{definition}

\theoremstyle{remark}
\newtheorem{rem}{Remark}

\newcommand{\entiere}[1]{\left\lfloor #1 \right\rfloor}

\newcommand{\abs}[1]{\left\lvert#1\right\rvert}

\setmarginsrb{2cm}{2cm}{2cm}{2.5cm}{0cm}{0cm}{0cm}{1cm}
\newcounter{ex}

\begin{document}

\date{}
\title{\textbf{A note about words which coincide except in one position}}
\author{Margot Bruneaux \\
{\small \texttt{margot.bruneaux@gmail.com}}}
\maketitle

\begin{abstract}
In this short note, we show that a result about words which coincide except in one position given as an exercise in Lothaire's \emph{Algebraic Combinatorics on Words} is false. Moreover, we derive a modified statement which allows us to fix the proof of a theorem which originally used the result of this exercise.
\end{abstract}

\bigskip

\bigskip

In \cite[p. 276]{Lot02}, the authors use in their proof of Theorem 8.1.11 the following statement.

\bigskip

\emph{Let $w=a_1a_2\cdots a_n$ and $v=b_1b_2\cdots b_n$ be two words having the same length $n$, such that $w$ has period $p$  and $v$ has period $q$ with $p\neq q$ and $p+q \leqslant n$. }

\emph{Suppose that there exists a position $t$, $1 \leqslant t \leqslant n$, such that for any $i\neq t$, $1 \leqslant i \leqslant n$, one has that $a_i=b_i$ (i.e., the two words $w$ and $v$ coincide  except, maybe, in position $t$).}

\emph{Then, $w$ and $v$ both have period $r=\mathrm{gcd}(p,q)$ and thus $w=v$.}

\bigskip

This result is not proved but only given as an exercise (\cite[Problem 8.1.4]{Lot02}). As it stands, this statement is false. Indeed, if one considers the words $w=ababab$ and $v=abaaab$ which coincide except in position 4 then $w$ has period $p=2$, $v$ has period $q=4$ and $p+q=6$ but $v$ has not period $\text{gcd}(2,4)=2$ and $w\neq v$. However, it becomes true by adding the bound $\lfloor \frac{n}{2} \rfloor$ for both periods as we can see with the following Proposition.

\bigskip

\begin{prop} Let $w$ and $v$ be two words having the same length $n$ such that $w$ has period $q$  and $v$ has period $p$. Assume that $w$ and $v$ coincide except, maybe, in one position. If $\max\{p,q\}\leqslant \left\lfloor\frac{n}{2}\right\rfloor$ then $w=v$. 
\end{prop}

\emph{Proof}. Without loose of generality, we can assume $q\leqslant p$. Moreover, if $p=q$ the result is clear so we can assume $q<p$.

Since $p\leqslant \left\lfloor\frac{n}{2}\right\rfloor$, by taking the reverse words of $w$ and $v$ if necessary, we can write 
\[w=a_1\cdots a_p\cdots x a_{r+1}\cdots a_{m} \qquad \text{and} \qquad v=a_1\cdots a_p\cdots y a_{r+1}\cdots a_{m}\]
with, maybe, $a_{r+1}\cdots a_m=\varepsilon$ is the empty word.

Let $c=\text{pgcd}(p, q)$. 

If $q=c$, the result is clear since the two words have then period $p$.

Assume $q \neq c$. Let $k$ and $k'$ be two integers such that $k'p+kq=-c$ with $k'<0$ and $k>0$.

We will prove that $y=a_{r+1-c}$.

We see $k$ and $k'$ as ``stockpiles of moves'' of lengths $q$ and $p$ taking into account the sign. For example, if $p=5$ et $q=3$, one can write $-2\times 5 + 3\times 3=-1$. Then, one has a stockpile of $3$ moves of length $3$ to the right and a stockpile of $2$ moves of length $5$ to the left.

One proceed the following way: starting from $y$, one makes alternately moves of length $p$ to the left in the word $v$ and moves of length $q$ to the right in the word $w$. Thus, for the above example, one obtains with $n=10$ and $r=5$:

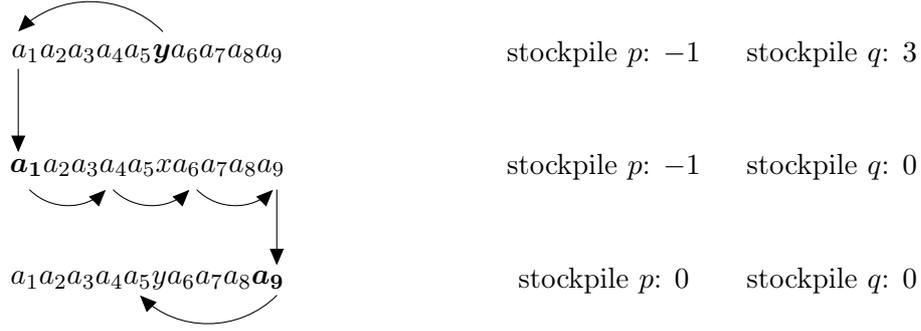
\begin{tikzpicture}[line cap=round,line join=round,>=triangle 45,x=1cm,y=1cm]
\clip(1,0.5) rectangle (20,6.5);
\draw[color=black] (5,5) node {$a_1a_2a_3a_4a_5\boldsymbol{y}a_6a_7a_8a_9$};
\draw[->] (5.2,5.3) to[bend left=-45] (3.3,5.3);
\draw[color=black] (5,3.5) node {$\boldsymbol{a_1}a_2a_3a_4a_5xa_6a_7a_8a_9$};
\draw[->] (3.3,4.8)-- (3.3,3.7);
\draw[->] (3.45,3.2) to[bend left=-45] (4.45,3.2);
\draw[->] (4.55,3.2) to[bend left=-45] (5.55,3.2);
\draw[->] (5.65,3.2) to[bend left=-45] (6.65,3.2);
\draw[color=black] (5,2) node {$a_1a_2a_3a_4a_5ya_6a_7a_8\boldsymbol{a_9}$};
\draw[->] (6.7,3.2)-- (6.7,2.2);
\draw[->] (6.7,1.8) to[bend left=45] (4.9,1.8);
\draw[color=black] (11,5) node {stockpile $p$: $-1$}; 
\draw[color=black] (14,5) node {stockpile $q$: $3$};
\draw[color=black] (11,3.5) node {stockpile $p$: $-1$};
\draw[color=black] (14,3.5) node {stockpile $q$: $0$};
\draw[color=black] (11,2) node {stockpile $p$: $0$};
\draw[color=black] (14,2) node {stockpile $q$: $0$};
\end{tikzpicture}

In general case, one proceeds in the same way: one makes as many moves of length $p$ as possible to the left in the word $v$ (at least one move is possible since $y$ stands after $a_p$). 

Then, one turns to the word $w$ and one makes as many moves of length $q$ as possible to the right  (at least one move is possible since $q<p$).

\begin{enumerate}
\item[\tiny$\bullet$] If one arrives at $w_{r+1}=x$ then $w_{r+1}=y$  and thus $x=y$.

By the Fine and Wilf theorem, since $n\geqslant p+q$, $c$ is a period of $w$ because $w=v$ and so $y=a_{r+1-c}$.
\item[\tiny$\bullet$] If one exhausts the stockpile of $q$ and arrives at a final position which is different from that of $x$, one returns to the word $v$ and makes all the remaining moves of length $p$ to the left. Thus, one arrives at $v_{r+1-c}=a_{r+1-c}$ and so $y=a_{r+1-c}$.
\item[\tiny$\bullet$] If one arrives at a different position from that of $x$ without having exhausted the stockpile of $q$, one returns to $v$. One makes as many moves of length $p$ as possible to the left (at least one move is possible since $p+q\leqslant n$). One arrives at a position which is strictly smaller than $r+1$ otherwise the final position will be strictly greater than $r+1-c$ which is impossible.

Then, one restarts the process. It will come to an end since, at each step, stockpiles of $p$ and $q$ strictly decrease in absolute terms.
\end{enumerate}

So $a_{r+1-c}=y$.

We reason in the same way to derive that $x=a_{r+1-c}$ by taking in this case B\'ezout's identity in the form $hp+h'q=-c$ with $h>0$ and $h'<0$. 

Thus, $y=a_{r+1-c}=x$ and so $w=v$. \hfill {\tiny{$\blacksquare$}}

\begin{rem} The above example $w=ababab$ and $v=abaaba$ shows the bound $\lfloor \frac{n}{2} \rfloor$ is the best possible.
\end{rem}

Theorem 8.1.11 of \cite{Lot02} is the following statement where $\Pi(w)$ denotes the set of all periods of a word $w$, with $0$ included.

\begin{thm} Let $\Pi=\{0=p_0 < p_1 <  \cdots < p_s=n\}$ be a set of integers and let $\delta_h=p_h-p_{h-1}$, $1 \leqslant h \leqslant s$. Then the following conditions are equivalent.
\begin{enumerate}
\item[(i)] There exists a word $w$ over  a two-letter alphabet with $\Pi(w)=\Pi$.
\item[(ii)] There exists a word $w$ with $\Pi(w)=\Pi$.
\item[(iii)] For each $h$, such that $\delta_h \leqslant n-p_h$, one has
\begin{enumerate}
\item[(a)] $p_h+k\delta_h \in \Pi$, for $k=1, ..., \lfloor (n-p_h)/\delta_h\rfloor$, and
\item[(b)] if $\delta_{h+1} < \delta_h$, then $\delta_h+\delta_{h+1}>n-p_h+\mathrm{gcd}(\delta_h,\delta_{h+1})$.
\end{enumerate}
\item[(iv)] For each $h$, such that $\delta_h \leqslant n-p_h$, one has
\begin{enumerate}
\item[(a)] $p_h+\delta_h \in\Pi$ and
\item[(b)] if $\delta_h=k\delta_{h+1}$, for some integer $k$ then $k=1$.
\end{enumerate}
\end{enumerate}
\end{thm}

In their proof, the authors refer to the result of Problem 8.1.4 to derive $(iv)$ implies $(i)$. More precisely, denoting $\Pi_h=\{p-p_h \mid p\in\Pi \text{ and } p\geqslant p_h\}$, $0\leqslant h \leqslant s$, they construct binary strings $w_h$ such that $\Pi(w_h)=\Pi_h$. For it, they use Problem 8.1.4 to prove the following result.

\begin{lem} If $\delta_h>n-p_h$ then there exists a sequence $a_1, ..., a_{\delta_h-\abs{w_h}}$ of letters in the same binary alphabet as $w_h$ such that the word $w_{h-1}=w_h a_1\cdots a_{\delta_h-\abs{w_h}} w_h$ has no period of length smaller than $\delta_h$.
\end{lem}

We will now see that Proposition 1 allows us to derive this lemma.

\bigskip

\emph{Proof of Lemma 1}. One uses mathematical induction on $m=\delta_h-\abs{w_h}$. 

Suppose $m=1$. Consider the two words $w_hxw_h$ and $w_hyw_h$ where $x$ and $y$ denote the two different letters in the binary alphabet of $w_h$. Assume by contradiction $w_hxw_h$ has period $p$ and $w_hyw_h$ has period $q$ such that $\max\{p, q\} < \delta_h$. Then, since $\delta_h=\abs{w_hx}$, $\max\{p, q\} \leqslant \lfloor \frac{\abs{w_hxw_h}}{2}\rfloor$ so, by Proposition 1, $x=y$ which is impossible since $x\neq y$. 

Assume the property true for a certain integer $m\geqslant 1$, i.e., there are letters $a_1$, ..., $a_m$ such that $w_ha_1\cdots a_mw_h$ has no period smaller than $\abs{w_h}+m$. Suppose that putting a letter $x$ or a letter $y$ (with $x\neq y$) between $a_{\lceil m/2 \rceil}$ and $a_{\lceil m/2 \rceil+1}$, we get two words that each have period smaller than or equal to $\abs{w_h}+\lceil m/2 \rceil$. Then, since $x\neq y$, one of these periods is strictly smaller than $\abs{w_h}+\lceil m/2 \rceil$.

The length of the two words we obtain is $L=2\abs{w_h}+m+1$ so $\entiere{L/2}=\abs{w_h}+\lceil m/2 \rceil$. Thus, by Proposition 1, $x=y$ which is absurd. 

So, there is a letter $b$ such that the word $w_{h-1}=w_ha_1\cdots a_{\lceil m/2 \rceil}b a_{\lceil m/2 \rceil+1}\cdots a_{m}w_h$ has no period smaller than or equal to $\abs{w_h}+\lceil m/2 \rceil$. Moreover, $2(\abs{w_h}+\lceil m/2 \rceil+1)>L$ then, by induction hypothesis, $w_{h-1}$ has no period smaller than $\delta_h=\abs{w_h}+m+1$. \hfill {\tiny{$\blacksquare$}}

\bibliographystyle{amsalpha}
\bibliography{biblio}  

\end{document}